\documentclass[12pt]{article}

\makeatletter
\@addtoreset{equation}{section}
\makeatother

\input epsf

\usepackage{amssymb}
\usepackage{amsmath}

\def\tfrac#1#2{{{\lower.6ex
\hbox{$\scriptstyle#1$}}\over
{\raise.7ex
\hbox{$\scriptstyle#2$}}}}

\def\erfc{{\rm erfc}}

\def\bigO{{\cal O}}

\def\protectbold#1{\protect{\boldmath{$#1$}}}
 
\def\Frac#1#2{\frac
{
 {\raise.6ex
 \hbox{$\disp#1$}}
}
{
 {\lower.6ex
 \hbox{$\disp#2$}}
 }
}

\begin{document}
\title{The leaky aquifer function revisited}
\author{
    Nico M. Temme\\
    CWI, Science Park 123, 1098 XG Amsterdam, The Netherlands. \\
     { \small e-mail: {\tt
    Nico.Temme@cwi.nl}}}
\date{}
\maketitle
\begin{center}
{\it Dedicated to Frank E. Harris on the occasion of his 80th birthday}
\end{center}

\begin{abstract}
This papers discusses  the leaky aquifer function considered in a recent paper by Frank Harris in the Journal of Computational and Applied Mathematics (2008). We describe properties of an integral representing this function and give details on how to compute this function with a single algorithm for a wide range of the parameters.
\end{abstract}

\section{Introduction}\label{intro}
The function 
\begin{equation}\label{intro1}
K_{\nu}(x,y)=\int_1^\infty e^{-xt-y/t}\frac{dt}{t^{\nu+1}}
\end{equation}
is related with the modified Bessel function (the MacDonald function)
\begin{equation}\label{intro2}
K_{\nu}(z)=\frac12\int_0^\infty e^{-z(t+1/t)}\frac{dt}{t^{\nu+1}}.
\end{equation}
For $\nu=0$ hydrologists call the function in \eqref{intro1} the {\it leaky aquifer function}, because this function can be used to describe water levels in pumped aquifer systems with finite transmissivity
and leakage could be analyzed in terms of this integral. 

In a recent paper by Frank Harris \cite{Harris:2008:IBG} the function \eqref{intro1} with general $\nu$ is considered and we refer to this paper for references to the hydrological literature. Also for $\nu=1,2,3,\ldots$ this function is useful in other hydrological systems.

The function in \eqref{intro1} can also be viewed as a generalization of the incomplete gamma function. In fact this is the starting point used in \cite{Chaudhry:1996:ACF} with the notation
\begin{equation}\label{intro3}
\Gamma(\alpha,x;b)=\int_x^\infty t^{\alpha-1} e^{-t-b/t}\,dt,
\end{equation}
which for $b=0$ is the standard incomplete gamma function
\begin{equation}\label{intro4}
\Gamma(\alpha,x)=\int_x^\infty t^{\alpha-1} e^{-t}\,dt,
\end{equation}
and it is easily verified that
\begin{equation}\label{intro5}
K_\nu(x,y)=x^\nu\Gamma(-\nu,x;xy).
\end{equation}
In \cite{Chaudhry:1996:ACF} closed forms are given in terms of modified Bessel functions and error functions when $\nu=\frac12-n$, $n=0,1,2,\ldots$.

We also mention a symmetry relation. We use in \eqref{intro1} the transformation $t\to 1/t$, which gives
\begin{equation}\label{intro6}
K_{\nu}(x,y)=\int_0^1 t^{\nu-1} e^{-x/t-yt}\,dt,
\end{equation}
from which we deduce (see also \eqref{intro2})
\begin{equation}\label{intro7}
K_{\nu}(x,y)+K_{-\nu}(y,x)=\int_0^\infty e^{-xt-y/t}\frac{dt}{t^{\nu+1}}=
2(x/y)^{\nu/2}K_\nu\left(2\sqrt{xy}\right).
\end{equation}

Frank Harris observed that the difference in the names assigned to the functions $K _\nu(x, y)$ and  $\Gamma(\alpha,x;b)$ was probably the reason there have been no previous communications connecting the several research communities. As application areas he mentions heat conduction, probability theory, electronic structure in periodic systems, and hydrology.

In his paper \cite{Harris:2008:IBG}  Harris gives a number of analytical properties of  $K _\nu(x, y)$ in the form of relations with other special functions, but the main part of the paper is an overview of expansions for this function. It presents several new expansions that are in various parameter ranges computationally more efficient than any of earlier proposed methods of evaluation. Several expansions are used in numerical evaluations for certain parameter ranges.

In particular we mention the small $y$ expansion in terms of incomplete gamma function
\begin{equation}\label{intro8}
K_{\nu}(x,y)=x^\nu\sum_{j=0}^\infty \Gamma(-\nu-j,x)\frac{(-xy)^j}{j!}.
\end{equation}
When $x$ is small it is better to use a similar expansion for $K_{-\nu}(y,x)$ combined with the relation in \eqref{intro7}, which requires the evaluation of the modified Bessel function.

In the present paper we exclude these small $x,y$ cases, assuming for example $x\ge1$ and $y\ge1$, and for other values we describe one numerical method based on numerical quadrature. We use a few transformations to bring the integral in a certain standard form, after which simple quadrature methods can be used for numerical evaluations. This method is not comparable in efficiency as some expansions derived in \cite{Harris:2008:IBG} for certain parameter ranges, but the main advantage is a universal method, that can be described in terms of a simple algorithm that can be used for very large parameter ranges of the function.

\section{Major contributions and scaling }\label{major}
First we like to see which part of the interval in \eqref{intro1} gives the main contribution. The maximal value of the exponential function occurs at $t=\sqrt{y/x}$. When $y>x$ this point is inside the interval of integration, and when $y<x$ it is outside this interval. In the latter case the exponential function is maximal at the  endpoint $t=1$ of the interval of integration.

When we include the parameter $\nu$ in this analysis, we write
\begin{equation}\label{major1}
K_{\nu}(x,y)=\int_1^\infty e^{\phi(t)}\frac{dt}{t},
\end{equation}
where
\begin{equation}\label{major2}
\phi(t)=-xt-y/t-\nu\ln t,
\end{equation}
with derivative
\begin{equation}\label{major3}
\phi^{\prime}(t)=-\frac{xt^2+\nu t-y}{t^2}.
\end{equation}
The function $\phi(t)$  is maximal at the zeros of $\phi^{\prime}(t)$. The relevant (positive) zero is
\begin{equation}\label{major4}
t_0=t_0(x,y,\nu)=\frac{\sqrt{\nu^2+4xy}-\nu}{2x}.
\end{equation}
We see from this form, or easier from \eqref{major3}, that $t_0=1$ when $y=x+\nu$. Hence, when $y>x+\nu$ the maximal value of $\phi(t)$ is inside the interval of integration at $t=t_0$, and when $y<x+\nu$ this function assumes its maximal value at the endpoint $t=1$. In particular when the parameters are large we conclude that a rough estimate of $K_\nu(x,y)$ is given by
\begin{equation}\label{major5}
K_{\nu}(x,y)=\left\{
\begin{array}{ll}
\bigO\left(e^{\phi(t_0)}\right) & \mbox {if \ } y>x+\nu \\
\bigO\left(e^{-x-y}\right)         & \mbox{if \ }  y \le x+\nu.
\end{array}
\right.
\end{equation}

In the $\{x,y,\nu\}$ parameter space, crossing the ``transition plane" $y=x+\nu$ causes a change in behavior of $K_{\nu}(x,y)$. Especially when the parameters are large, the function changes from small values $\exp(\phi(t_0))$ at one side of the plane to very small values $\exp(-x-y)$ at the other side of the plane.

Considering this from the viewpoint of asymptotic analysis, we call $t_0$ a saddle point, and when this point is properly inside the interval of integration we can use Laplace's method (see \cite[\S II.1]{Wong:2001:AAI}) to obtain asymptotic representations of which the dominant factor is given in the first case of \eqref{major5}. When the saddle point $t_0$ is properly inside $[0,1]$ we can use Watson's lemma (see \cite[\S I.5]{Wong:2001:AAI}), after a transformation to a Laplace-type integral. Although the asymptotic forms in \eqref{major5} pass continuously into each other, the detailed asymptotic forms obtained after applying Laplace's method or Watson's lemma are discontinuous at the transition plane $y=x+\nu$.

In \cite{Chaudhry:1996:ACF} we have given an asymptotic representation of $\Gamma(\alpha,x;b)$  in terms of the complementary error function
\begin{equation}\label{major6}
\erfc\,z=\frac{2}{\sqrt{\pi}}\int_z^\infty e^{-t^2}\,dt,
\end{equation}
which representation behaves smoothly at the transition manifold (in that case defined by the relation $b=x(x-\alpha)$), and which is valid in a wide parameter domain of the three variables $\alpha, x, b$. The same result  can be written in terms of the present parameters $x,y,\nu$ of the function $K_\nu(x,y)$. 

Understanding the global behavior as given in \eqref{major5} of the function to be computed is important when choosing a numerical algorithm, because, in particular when the parameters $x,y,\nu$ are large, the function values may be extremely small, and numerical instabilities may arise in certain representations. For example, scaling the function before a numerical quadrature method is used is important when relative accuracy is wanted. In that case we can write, using \eqref{major1} and \eqref{major5},
\begin{equation}\label{major7}
K_{\nu}(x,y)=e^{\phi(t_m)}\int_1^\infty e^{\phi(t)-\phi(t_m)}\frac{dt}{t},
\end{equation}
where $t_m=t_0$ (given in \eqref{major4}) when $y>x+\nu$ and $t_m=1$ when $y\le x+\nu$.

In a similar way, if we prefer the representation in \eqref{intro6} as the starting point for numerical evaluations (as we do), we can write 
\begin{equation}\label{major8}
K_{\nu}(x,y)=e^{\psi(t_m)}\int_0^1 e^{\psi(t)-\psi(t_m)}\frac{dt}{t},\quad
\psi(t)=-x/t-yt+\nu\ln t,
\end{equation}
where $t_m=t_0(y,x,-\nu)$ (with $t_0$ given in \eqref{major4}) when $y>x+\nu$ and $t_m=1$ when $y\le x+\nu$.

\section{Numerical quadrature of \protectbold{K_\nu(x,y)}}\label{quad}
In the numerical literature  the role of the trapezoidal rule has been discussed with a number of examples, from which it is concluded that this rule gives excellent performances for a certain type of integrals, and how given integrals can be transformed into this type. In particular, when the functions to be integrated are analytic (as in the present case) several efficient transformations are discussed. Detailed information, also with a discussion on error estimations and  exact error representations in the form of contour integrals, can be found in the work by Mori and co-workers 
\cite{Mori:2001:DET,Mori:2005:DDE,Mori:1973:QVT,Takahasi:1970:EEN,Mori:1974:DEE}, and in the recent books \cite[\S4.2]{Kythe:2005:HCM} and \cite[\S 5.4]{Gil:2007:NSF}.

The trapezoidal rule works with great efficiency for integrals of the type $\int_0^1 f(t)\,dt$ when $f$ is very smooth, say analytic, and if  for the derivatives hold $f^{(n)}(0)=f^{(n)}(1)$. In that case all terms in the error estimate of this quadrature rule based on the Euler-Maclaurin formula (see \cite[p.~131]{Gil:2007:NSF}) vanish, and we have convergence of exponential type.

In the present integral in \eqref{major8} we have vanishing derivatives of all orders at $t=0$, but at $t=1$ this property does not hold. So, straightforward  application of the compound trapezoidal rule to \eqref{major8} will  give the standard  error estimate $\bigO(h^2)$, where $h$ is the stepsize of the rule.

We use the following simple substitution
\begin{equation}\label{quad1}
t=\tanh \frac {s}{1-s},\quad 0\le s \le 1,
\end{equation}
which gives 
\begin{equation}\label{quad2}
K_{\nu}(x,y)=e^{\psi(t_m)}\int_0^1 \frac{e^{\psi(t)-\psi(t_m)}}{t\,(\cosh(s/(1-s))\,(1-s))^2}\,ds.
\end{equation}
In this representation all derivatives at both endpoints vanish. For small values of $x$ the effect of vanishing derivatives at $s=0$ is less noticeable, and we refer in that case to the efficient expansions mentioned at the end of \S\ref{intro}. Small values of $y$ are allowed, however.  

It is possible to use different substitutions, but we continue with \eqref{quad1}. For other transformations of the integral as in \eqref{major8} we refer to the papers by Mori and co-workers, with one based on the error function \eqref{major6}. That rule requires in each function evaluation the evaluation of the error function, but it may have excellent convergence properties.

\renewcommand{\arraystretch}{1.0}
\begin{table}
\caption{Function values of $K_\nu(x,y)$ and relative errors  $\delta$ for a selection of parameter values $x,y,\nu$. In the first  nine rows $h=1/40$, in the last three rows $h=1/80$.
\label{table1}}
$$
\begin{array}{rrrrrr}\hline
x \ & y \ &  \nu \  & K_\nu(x,y) \quad\quad\quad& \delta  \ \quad\quad \\
\hline\\[-10pt]
4.95 & 5.00 & 2.00 & 0.12249\,98798\times10^{-004} & 0.16\times10^{-10} \\
10.0 & 2.00 & 6.00 & 0.41500\,45943\times10^{-006} & 0.63\times10^{-10} \\
3.10 & 2.60 & 5.00 & 0.52850\,43253\times10^{-003} & 0.27\times10^{-10} \\
&&&&\\
49.0 & 50.0 & 20.0 & 0.44311\,56799\times10^{-044} & 0.96\times10^{-10}\\
100.0 & 20.0 & 60.0 & 0.54438\,05280\times10^{-054} & 0.60\times10^{-09} \\
31.0 & 26.0 & 50.0 & 0.31405\,73138\times10^{-026} & 0.24\times10^{-09}\\
&&&&\\
490.0 & 500.0 & 200.0 & 0.57348\,63502\times10^{-432} & 0.56\times10^{-10}\\
1000.0 & 200.0 & 600.0 & 0.5014\,041537\times10^{-524} & 0.92\times10^{-04}\\
310.0 & 260.0 & 500.0 & 0.51400\,54359\times10^{-250} & 0.60\times10^{-06}\\
&&&&\\
490.0 & 500.0 & 200.0 & 0.57348\,63503\times10^{-432} & 0.36\times10^{-12}\\
1000.0 & 200.0 & 600.0 & 0.50145\,04964\times10^{-524} & 0.47\times10^{-12}\\
310.0 & 260.0 & 500.0 & 0.51400\,57464\times10^{-250} & 0.75\times10^{-12}\\ \hline
\end{array}
$$
\end{table}
\renewcommand{\arraystretch}{1.0}

To give a an idea about this starting point for applying the trapezoidal rule, we use Cases 2, 3, and 4 considered in \cite{Harris:2008:IBG} and larger values of $x,y,\nu$.  Case 1 has $x=0.01$, which value we have excluded. In the first  9 rows of the Table~\ref{table1} we use 40 integrand evaluations; this means, we use the compound trapezoidal rule with $h=1/40$. This gives rather uniform relative errors $\delta$, except in rows 7-- 9. In rows 10 --12 we  use $h=1/80$ and we notice that again a rather uniform relative error $\delta$ is obtained. All computations are done with Maple, with Digits =15 (working precision of 15 significant digits).

\subsection{About error estimations}\label{errors}
As mentioned earlier, Mori and co-authors have given estimates and representations of the error in the application of the trapezoidal rule for integrals of the form \eqref{quad2}. Let us write the integral in the form
\begin{equation}\label{quad3}
I=\int_0^1f(s)\,ds.
\end{equation}
Then one estimate of the error $E_h$ in the trapezoidal rule
\begin{equation}\label{quad4}
I=T_h+E_h, \quad T_h=\tfrac12h(f(0)+f(1))+h\sum_{k=0}^{n-1}f(kh), \quad h=1/n,
\end{equation}
can be written in the form of 
\begin{equation}\label{quad5}
E_h \sim -\int_{{\cal C}_+} f(s) e^{2\pi i s/h}\,ds-\int_{{\cal C}_-} f(s) e^{-2\pi i s/h}\,ds,
\end{equation}
where ${\cal C}_+$ is a contour in the half-plane $\Im s>0$ and ${\cal C}_-$ in $\Im s<0$, both running from $0$ to $1$. The contours may be close to the original interval $[0,1]$, but they can be deformed in order to pick up the relevant information from saddle points of the integrand or singular points of the function $f(s)$.

Similar estimates follow from Fourier representations of $E_h$. We have \cite[\S15.2]{Luke:1969:SFB}
\begin{equation}\label{quad6}
E_h  =-2 \sum_{k=1}^\infty F_k,\quad F_k=\int_0^1f(s) \cos(2\pi ks/h)\,ds.
\end{equation}
Because $f(s)$ has vanishing derivatives of all orders, this Fourier series is rapidly convergent and $E_h$ can be approximated by the first term of the series. That is, $E_h\sim-2F_1$. By writing the cosine as exponentials, and deforming the interval $[0,1]$ into contours in the complex plane, we obtain the estimate in \eqref{quad5}.

For simple functions it is possible to obtain estimates of the integrals in \eqref{quad5} for small values of $h$ from saddle points or singular points. We are dealing, however, with a function that depends on three parameters $x, y, \nu$, which may be large as well (we mean, much larger than $1/h$). In addition, we use a transformation from the integral in \eqref{major8} into \eqref{quad1} by \eqref{quad1}, which is a rather simple transformation, but for a saddle point analysis it causes extra complications.

A detailed saddle point analysis falls outside the scope of this paper. However, we can use a numerical method for estimating the error. An aspect of efficiency of the trapezoidal rule is the fact that we can use earlier computed function values when we halve the stepsize. After computing $T_h$ of \eqref{quad4} we use the shifted rule
\begin{equation}\label{quad7}
I= S_h+\widetilde E_h,\quad S_h=h\sum_{k=0}^{n-1} f\left(kh+\tfrac12h)\right), \quad \widetilde E_h=-2\sum_{k=1}^\infty (-1)^kF_k,
\end{equation}
and compute the absolute error $|T_h-S_h|$ (or the relative error). When this quantity satisfies a certain precision criterion, we stop the  computation by accepting the result $I\doteq T_{h/2}=\frac12(T_h+S_h)$. Because this is the trapezoidal rule with $h/2$ the value $T_{h/2}$ is a much better approximation in a fast converging rule than each of $T_h$ and $S_h$. Notice that $T_h-S_h\doteq-4F_2$, which is two times the approximate error in the rule for $T_{h/2}$.

We see this effect in Table~\ref{table1} in the two evaluations for $x=1000.0$, $y= 200.0$, $\nu=600.0$. For $h=1/40$ the relative error is   $0.92\times10^{-4}$, when we halve the stepsize it is $0.47\times 10^{-12}$.

We conclude that  if, for a certain $h$, we have $|T_h-S_h|<\epsilon$ we can be rather sure that for the rule for $T_{h/2}$ we have $|E_{h/2}|<\epsilon$; similar for relative errors.

\section{Concluding remarks}\label{concrem}
We have described a rather uniform method for computing the leaky aquifer function $K_\nu(x,y)$ based on a simple quadrature rule. As shown in detail in Frank Harris's paper \cite{Harris:2008:IBG},  by using expansions of this function that are designed for specific parameter ranges more efficient algorithms are possible than the one considered in the present paper.  

It is also true that this function, given in the form of an integral with positive integrand, is not difficult to compute by using standard quadrature routines as, for example, provided by packages as Matlab, Mathematica, or Maple, or by using routines written in Fortran in standard quadrature libraries. 

However,  by preparing the integral by using simple transformations,  an efficient quadrature rule can be used, which can easily be written as an algorithm and can be used  for high performance computing for a wide range of the parameters.

\bigskip
\noindent \textbf{Acknowledgments.\ {}} The author thanks the referee for valuable comments on the first version of the paper.
The author acknowledges
financial support of the Spanish \textit{Ministerio de Educaci\'{o}n
y Ciencia}, project MTM2006--09050, and of the \textit{Gobierno de
Navarra}, Res. 07/05/2008.

\end{document}